\documentclass{article}

\usepackage{arxiv}

\usepackage[utf8]{inputenc} 
\usepackage[T1]{fontenc}    
\usepackage{hyperref}       
\usepackage{url}            
\usepackage{booktabs}       
\usepackage{amsfonts}       
\usepackage{nicefrac}       
\usepackage{lipsum}
\usepackage{enumerate}

\usepackage{color}
\usepackage{multirow}
\usepackage{latexsym}
\usepackage{amsmath,amssymb,amsthm,graphicx}
\usepackage{dsfont}
\usepackage{enumitem}
\usepackage{hyperref}
\usepackage{comment}

\newtheoremstyle{thm}
{9pt}
{9pt}
{\itshape}
{}
{\bfseries}
{.}
{ }
{}
\theoremstyle{thm}

\newtheorem{theorem}{Theorem}[section]
\newtheorem{lemma}[theorem]{Lemma}
\newtheorem{corollary}[theorem]{Corollary}

\theoremstyle{definition}
\newtheorem{remark}[theorem]{Remark}

\newcommand{\vertk}{\stackrel{{\cal D}}{\longrightarrow}}

\theoremstyle{definition}

\newtheoremstyle{def}
{9pt}
{9pt}
{}
{}
{\bfseries}
{.}
{ }
{}
\theoremstyle{def}

\newcommand{\RR}{\mathbb{R}} 
\newcommand{\NN}{\mathbb{N}} 
\newcommand{\E}{\mathbb{E}} 
\newcommand{\PP}{\mathbb{P}} 

\newcommand{\jj}{\mathbf{j}}
\newcommand{\VV}{\mathcal{V}_n}
\newcommand{\PPP}[1]{\operatorname{\mathbb{P}}\left(\,#1\,\right)}

\newcommand{\ind}{1\hspace{-0,9ex}1}

\renewcommand{\footnoterule}{%
	\kern -3.5pt
	\hrule width \textwidth height 1pt
	\kern 3.5pt
}

\makeatletter
\def\blfootnote{\xdef\@thefnmark{}\@footnotetext}
\makeatother

\title{Limit laws for large $k$th-nearest neighbor balls}

\author{
Nicolas Chenavier\\
Universit\'e du Littoral C\^ote d'Opale\\
50 rue F. Buisson 62228 Calais\\
\texttt{nicolas.chenavier@univ-littoral.fr}\\
\And
Norbert Henze\\
Institute of Stochastics \\
Karlsruhe Institute of Technology (KIT) \\
Englerstr. 2, D-76133 Karlsruhe \\
\texttt{Norbert.Henze@kit.edu}\
\And  Moritz Otto\\
Institute for Mathematical Stochastics \\
Otto von Guericke University Magdeburg \\
Universit\"atsplatz 2, D-39106 Magdeburg \\
\texttt{Moritz.Otto@ovgu.de}\\
}

\begin{document}

\date{\today}
\maketitle


\begin{abstract}
Let $X_1,X_2, \ldots, X_n$ be {a sequence of independent random} points in $\mathbb{R}^d$ with common Lebesgue density $f$.
Under some conditions on $f$, we obtain a Poisson limit theorem, as $n \to \infty$, for the number of large probability $k$th-nearest
neighbor balls of $X_1,\ldots,X_n$. Our result generalizes Theorem 2.2 of \cite{GHW19}, which refers to the special case $k=1$. Our proof is completely different since it employs the Chen-Stein method instead
of the method of moments. Moreover, we obtain a rate of convergence for the Poisson approximation.
\end{abstract}

\section{Introduction and main results}
The starting point of this paper is the following result, see \cite{GHW19}. Let $X,X_1,\ldots,X_n, \ldots $ be a sequence of independent and identically distributed
(i.i.d.) random points in $\RR^d$ that are defined on a common probability space $(\Omega,{\cal A},\mathbb{P})$.
We assume that the distribution of $X$, which is denoted by $\mu$, is absolutely continuous with respect to Lebesgue measure $\lambda$,
and we denote the density of $\mu$ by $f$. Writing $\|\cdot \|$ for the Euclidean norm in $\RR^d$, and putting ${\cal X}_n := \{X_1,\ldots,X_n\}$, let $R_{i,n}:=\min_{j\ne i, j\le n} \|X_i-X_j\|$ be the distance
from $X_i$ to its nearest neighbor in the set ${\cal X}_n \setminus \{X_i\}$. Moreover, let $\ind \{A\}$ denote the indicator function of a set $A$, and write
$B(x,r) = \{y \in \RR^d: \|x-y\|\le r\}$ for the closed ball centered at $x$ with radius $r$. Finally, let
\[
C_n := \sum_{i=1}^n \ind \Big{\{} \mu \big( B(X_i,R_{i,n})\big) > \frac{t + \log n}{n} \Big{\}}
\]
denote the number of exceedances of probability volumes of nearest neighbor balls that are larger than the threshold $(t + \log n)/n$.
The main result of \cite{GHW19} is
Theorem 2.2 of that paper, which states that, under a weak condition on the density $f$, for each fixed    $t \in \RR$, we have
\begin{equation}\label{convnn}
C_n \vertk {\rm Po}(\exp(-t))
\end{equation}
as $n \to \infty$, where $\vertk$ denotes convergence in distribution, and Po$(\xi)$ is the Poisson distribution with parameter
$\xi >0$.

Since the {\em maximum} probability content of these nearest balls, denoted by $P_n$, is at most $(t+ \log n)/n$ if, and only if, $C_n =0$,
we immediately obtain a Gumbel limit $\lim_{n\to \infty} \mathbb{P}(nP_n - \log n \le t) = \exp(-\exp(-t))$ for $P_n$.

To state a sufficient condition on $f$ that guarantees \eqref{convnn}, let
$\text{supp}(\mu) := \{x \in \RR^d: \mu (B(x,r)) > 0 \textrm{ for each } r >0\}$ denote the support of $\mu$.
Theorem 2.2 of \cite{GHW19} requires that there are $\beta \in (0,1)$, $c_{max}<\infty$ and $\delta >0$
such that, for any $r,s >0$ and any $x,z \in$ supp$(\mu)$ with
$\|x-z\| \ge \max\{r,s\}$ and  $\mu\left( B(x,r)\right) = \mu\left( B(z,s)\right)\le \delta$, one has
\[
\frac{\mu\left( B(x,r)\cap   B(z,s)\right)}{\mu\left( B(z,s)\right)}
\le \beta
\]
and $\mu\big( B(z, 2s)\big)\le c_{max}\mu\left( B(z,s)\right)$.

These conditions hold if $\text{supp}(f)$ is a compact set $K$ (say), and there are
$f_-,f_+ \in (0,\infty)$ such that
\begin{equation}\label{densitybounds}
f_-\leq f(x)\leq f_+, \qquad x \in K.
\end{equation}
Thus, the density $f$ of $X$ is bounded and bounded away from zero.

The purpose of this paper is to generalize \eqref{convnn} to $k$th-nearest neighbors, and to derive a rate of convergence for the Poisson approximation of {the number of exceedances}. 

Before stating our main results, we give some more notation.  For fixed $k \le n-1$, we denote by $R_{i,n,k}$ the Euclidean
distance of  $X_i$ to its $k$th-nearest neighbor among ${\cal X}_n \setminus \{X_i\}$,
and we write  $B(X_i,R_{i,n,k})$ for the $k$th-nearest neighbor ball centered at $X_i$ with radius $R_{i,n,k}$.  For fixed $t \in \RR$, put
\begin{equation} \label{defynk}
v_{n,k} :=v_{n,k}(t) :=\frac{t+\log n+(k-1)\log\log n-\log(k-1)!}{n},
\end{equation}
and let
\begin{equation}\label{defcnk}
C_{n,k} := \sum_{i=1}^n \ind \big{\{}\mu \big( B(X_i,R_{i,n,k})\big) > v_{n,k}\big{\}}
\end{equation}
denote the {\textit{number of exceedances}} of probability contents of $k$th-nearest neighbor balls over the threshold $v_{n,k}$ defined in \eqref{defynk}.

The threshold $v_{n,k}$ is in some sense {\em universal} in dealing with the number of exceedances of probability contents of $k$th-nearest neighbor balls. To this end, suppose that, in
much more generality than considered so far, $X,X_1,X_2, \ldots $ are i.i.d. random elements taking values in a separable metric space $(S,\rho)$.
We retain the notations $\mu$ for the distribution of $X$ and $B(x,r) := \{y \in S: \rho(x,y) \le r\}$ for the closed ball with radius $r$ centered at $x \in S$. Regarding the distribution
$\mu$, we assume that
\begin{equation}\label{conditionbound}
\mu(\{y \in S: \rho(x,y) =r\}) =0, \qquad x \in S, \ r \ge 0.
\end{equation}
As a consequence, the distances $\rho(X_i,X_j)$, where $j \in \{1,\ldots,n\} \setminus \{i\}$, are different with probability    one
for each $i \in \{1,\ldots,n\}$. Thus, for fixed $k \le n-1$, there is almost surely a unique $k$th-nearest neighbor of $X_i$,
and we also retain the notations $R_{i,n,k}$ for the distance of  $X_i$ to its $k$th-nearest neighbor among ${\cal X}_n \setminus \{X_i\}$ and
$B(X_i,R_{i,n,k})$ for the ball centered at $X_i$ with radius $R_{i,n,k}$. Notice that the condition \eqref{conditionbound} excludes discrete metric spaces (see, e.g., Section 4 of \cite{ZO18}),
but not function spaces like, e.g., the space $C[0,1]$ of continuous functions on $[0,1]$ with the supremum metric, and with Wiener measure $\mu$. 

\medskip

{In what follows, for sequences $(a_n)_{n\geq 0}$ and $(b_n)_{n\geq 0}$ of real numbers, write $a_n=O(b_n)$ if  $|a_n|\leq C|b_n|$, $n\geq 1$, for some positive constant $C$.}

\begin{theorem}\label{thmmetricspace}
If $X_1,X_2, \ldots$ are i.i.d. random elements of a metric space $(S,\rho)$, and if \eqref{conditionbound} holds, then
the sequence $(C_{n,k})$ satisfies
\[
 \E[C_{n,k}] = {\rm e}^{-t} + O\left( \frac{\log\log n}{\log n}  \right).
\]
\end{theorem}

\medskip

In particular, the mean number of exceedances $C_{n,k}$ converges to ${\rm e}^{-t}$ as $n$ goes to infinity.
 By Markov's inequality, this result implies the tightness of the sequence $(C_{n,k})_{n \ge 1}$. Thus, at least a subsequence converges in distribution.
The next result states convergence of $C_{n,k}$ to a Poisson distribution if $(S,\rho) = (\RR^d,\|\cdot\|)$ and \eqref{densitybounds} holds.
To this end, let $d_{TV}(Y,Z)$ be the total variation between two integer-valued random variables $Y$ and $Z$, i.e.,
\[
d_{TV}(Y, Z) = 2\sup_{A\subset \NN}|\PP(Y\in A) -\PP(Z\in A)|.
\]

\begin{theorem}\label{Th:poissonapprox}
Let $Z$ be a Poisson random variable with parameter ${\rm e}^{-t}$. If $X, X_1,X_2, \ldots $ are i.i.d. in $\RR^d$ with density $f$, and if the distribution $\mu$ of $X$ has compact support {$[0,1]^d$}
 and satisfies \eqref{densitybounds}, then, as $n \to \infty$,
 \[
 d_{TV}\left(C_{n,k}, Z\right) =  O\left( \frac{\log\log n}{\log n} \right).
 \]
\end{theorem}
 Theorem \ref{Th:poissonapprox} is not only a generalization of Theorem 2.2 of \cite{GHW19} over all $k\geq 1$: it also provides a rate of convergence for the Poisson approximation of  $C_{n,k}$. {Our theorem is stated in the particular case {that} the support of $\mu$ is $[0,1]^d$ but we think that it can be extended to any measure $\mu$ whose support is a general convex body. For {the sake of} readibility of the manuscript, we did not deal with such a generalization.  }
 
 {\begin{remark}
The study of extremes of $k$th-nearest neighbor balls is classical in stochastic geometry, and it has various applications, see e.g. \cite{P}. In Section 4 in \cite{O}, bounds for the total variation distance of the process of Poisson points with large $k$th-nearest neighbor ball (with respect to the intensity measure) and a Poisson process were obtained. Parallel to our work, these results have been extended by  Bobrowski \textit{et al.} to the Kantorovich-Rubinstein distance and generalized to the binomial process in Section 6.2 of a paper  which has just been submitted \cite{BSY}. While the results in \cite{BSY}  and \cite{O} rely on Palm couplings of a thinned Poisson/binomial process and employ distances of point processes, we derive a bound on the total variation distance of the number of large $k$th-nearest-neighbor balls and a Poisson-distributed random variable. Our approach permits to build arguments  on classical Poisson approximation theory \cite{AGG} and  an asymptotic independence property stated in Lemma \ref{Le:independence} below, and it thus results in a considerably shorter and less technical proof.
\end{remark}} 

Now, let
 $P_{n,k} = \max_{1 \le i \le n} \mu\big(B(X_i,R_{i,n,k})\big)$ be the maximum probability content of the $k$th-nearest neighbor balls.
Since $C_{n,k} = 0$ if, and only if, $P_{n,k} \le v_{n,k}$, we obtain the following corollary.

\begin{corollary}
Under the conditions of Theorem~\ref{Th:poissonapprox}, we have
\[
\lim_{n\to \infty} \mathbb{P}\big(n P_{n,k} - \log n - (k-1)\log \log n + \log(k-1)! \le t\big) = {\rm G}(t), \quad t \in \mathbb{R},
\]
where ${\rm G}(t) = \exp(-\exp(-t))$ is the distribution function of the Gumbel distribution.
\end{corollary}

\begin{remark}
If, in the Euclidean case, the density $f$ is continuous, then  $\mu(B(X_i,R_{i,n,k}))$  is approximately equal to $f(X_i)\kappa_d R_{i,n,k}^d$, where
$\kappa_d = \pi^{d/2}/\Gamma(1+d/2)$ is the volume of the unit ball in $\RR^d$. Under additional smoothness assumptions on $f$ and  \eqref{densitybounds},
\cite{HE82, HE83} proved that
\begin{equation}\label{gumbelknnbound}
\lim_{n\to \infty} \PP\left( \max_{i=1,\ldots,n} f(X_i) \kappa_d \min\big(R_{i,n,k}^d, \|X_i-\partial K\|^d\big) \le v_{n,k}\right) = {\rm G}(t),
\end{equation}
\end{remark}
{where $K$ is the support of $\mu$}. Here, the distance $\|X_i-\partial K\|$ of $X_i$ to the boundary of $K$ is important to overcome edge effects. These effects dominate the asymptotic behavior
of the maximum of the $k$th-nearest neighbor distances if $k \ge d$, see  \cite{DH89, DH90}. In fact, \cite{HE82} proved convergence of the factorial moments of
\[
\widetilde{C}_{n,k} := \sum_{i=1}^n \ind \Big{\{} f(X_i)\kappa_d \min\big(R_{i,n,k}^d,\|X_i - \partial K\|^d\big) > v_{n,k} \Big{\}}
\]
to the corresponding factorial moments of a random variable with the Poisson distribution ${\rm Po}({\rm e}^{-t})$ and thus, by the method of moments, more than   \eqref{gumbelknnbound}, namely
$\widetilde{C}_{n,k} \overset{\mathcal{D}}{\rightarrow} {\rm Po}\big({\rm e}^{-t}\big)$. However, our proof of Theorem~\ref{Th:poissonapprox} is completely different thereof, since it is based on the
Chen-Stein method and provides a rate of convergence.

\section{Proofs}

\subsection{Proof of Theorem \ref{thmmetricspace}}
\begin{proof}
By symmetry, we have
\begin{eqnarray*}
\mathbb{E} \big{[}C_{n,k}\big{]} & =  & n\,  \mathbb{P} \left( \mu (B(X_1,R_{1,n,k})) > v_{n,k} \right)\\
& = & n \, \mathbb{E} \left[ \mathbb{P} \left( \mu(B(X_1,R_{1,n,k}))\right) > v_{n,k}|X_1  \right].
\end{eqnarray*}
For a fixed $x \in S$, let
\[
H_x(r):=\mathbb{P} \left(\rho(x,X) \le r\right), \quad r \ge 0,
\]
be the distribution function of $\rho(x,X)$. Due to the condition \eqref{conditionbound}, the function $H_x$ is continuous,
and by the probability integral transform (see e.g. \cite{BD15}, p. 8),
the random variable
\[
H_x(\rho(x,X)) = \mu \big( B(x,\rho(x,X))\big)
\]
is uniformly distributed in the unit interval $[0,1]$.
Put $U_j := H_x(\rho(x,X_{j+1}))$, $j=1,\ldots,n-1$. Then $U_1,\ldots,U_{n-1}$ are i.i.d. random variables with
a uniform distribution in $(0,1)$. Hence, conditionally on $X_1=x$, the random variable $\mu (B(X_1,R_{1,n,k}))$ has the same distribution as
$U_{k:n-1}$, where $U_{1:n-1} < \ldots < U_{n-1:n-1}$ are the order statistics of $U_1,\ldots,U_{n-1}$, and this distribution does not depend on $x$.
Now, because of a well-known relation between the distribution   of order statistics from the uniform distribution on $(0,1)$ and the binomial distribution (see, e.g., \cite{ANS}, p. 16),
we have
\[
\mathbb{P}(U_{k:n-1} > s) = \sum_{j=0}^{k-1} {n-1 \choose j} s^j (1-s)^{n-1-j}
\]
and thus
\begin{equation}\label{expcnk}
\mathbb{E} \big{[}C_{n,k}\big{]} = n \sum_{j=0}^{k-1} {n-1 \choose j} v_{n,k}^j (1-v_{n,k})^{n-1-j}.
\end{equation}
Here, the summand for $j=k-1$ equals
\[n {n-1 \choose k-1} v_{n,k}^{k-1} (1-v_{n,k})^{n-k} = \frac{n}{(k-1)!}\, (nv_{n,k})^{k-1}\, \prod_{i=1}^{k-1}\frac{n-i}{n}\, (1-v_{n,k})^{n-k}.\]
Using Taylor expansions, \eqref{defynk} yields
 \[nv_{n,k}  = \log n+O\left(\log\log n\right), \quad
 \prod_{i=1}^{k-1}\frac{n-i}{n} = 1+O\left(  \frac{1}{n}\right)
 \]
 and
\[
(1-v_{n,k})^{n-k} = \frac{(k-1)!}{n}\, \exp\left(-t-(k-1)\log\log n+O\left( \frac{\log^2(n)}{n}  \right)\right).
\]
Straigthforward computations now give
\[n {n-1 \choose k-1} v_{n,k}^{k-1} (1-v_{n,k})^{n-k}  = {\rm e}^{-t}+O\left( \frac{\log\log n}{\log n}  \right).\]
Regarding the remaining summands on the right hand side of \eqref{expcnk}, it is readily seen that
\[
\sum_{j=0}^{k-2} {n-1 \choose j} v_{n,k}^j (1-v_{n,k})^{n-1-j} = O\left( n {n-1 \choose k-1} v_{n,k}^{k-1} (1-v_{n,k})^{n-k} \cdot \frac{1}{nv_{n,k}}  \right),
\] with the convention that the sum equals 0 if $k=1$. From the above computations and from  \eqref{defynk}, it follows that this sum
equals $O\left(1/\log n \right)$, which concludes the proof of Theorem \ref{thmmetricspace}.
\end{proof}

\begin{remark}\label{rem-byproduct}
In the proof given above, we conditioned on the realizations $x$ of $X_1$. Since the distribution of $H_x(\rho(x,X)) = \mu\big(B(x,\rho(x,X))\big)$ does
not depend on $X$, we obtain as a by-product that
\[
\PP\left(\mu (B(X_1,R_{1,n,k})) > v_{n,k} \right)  = \sum_{j=0}^{k-1} {n-1 \choose j} v_{n,k}^j (1-v_{n,k})^{n-1-j} \sim \frac{{\rm e}^{-t}}{n},
\]
if $X_1,X_2,\ldots,X_n$ are independent and $X_2,\ldots,X_n$ are i.i.d. according to $\mu$. Here, $X_1$ may have an arbitrary distribution and $a_n \sim b_n$ means that $a_n/b_n \to 1$ as $n \to \infty$.
\end{remark}

\subsection{Proof of Theorem~\ref{Th:poissonapprox}}
The main idea to derive Theorem \ref{Th:poissonapprox} is to discretize {$\text{supp}(\mu)=[0,1]^d$} into finitely many ``small sets''
and then to employ the Chen-Stein method. To apply this method, we will have to check an  \textit{asymptotic independence property}
and a \textit{local property} which ensures that, with high probability,  two exceedances cannot appear in the same neighborhood.
We introduce these properties below and recall a result due to Arratia \textit{et al.} \cite{AGG} on the Chen-Stein method.

\paragraph{The asymptotic independence property}
Fix $\varepsilon > 0$. Writing $\lfloor \cdot \rfloor$ for the floor function, we partition {$[0,1]^d$} into a set $\VV$ of $N_n^d$  \textit{subcubes} (i.e., subsets that are cubes)
 of equal size that can only have boundary points in common, where $N_n=\lfloor n/\log(n)^{1+\varepsilon}\rfloor$. The subcubes are indexed
by the set  $\left[ 1,N_n\right] ^{d} = \{\jj:=(j_{1},\ldots, j_{d}): j_m \in \{1, \ldots, N_n\} \text{ for } m \in \{1,\ldots,d\}\}$.
With a slight abuse of notation, we identify a cube with its index.  Let
\[
\mathcal{E}_n = \bigcap_{\jj\in \VV}\{{\cal X}_n\cap \jj\neq \emptyset\}
\]
be the event that each of the subcubes contains at least one of the points of ${\cal X}_n$.
    The event $\mathcal{E}_n$ is extensively used in stochastic geometry to derive central limit theorems or to deal with extremes \cite{AB, BC, CR},
    and it will play a crucial role throughout  the rest of the paper. The following lemma, which captures the idea of ``asymptotic independence''  ,
     is at the heart of our development.

\begin{lemma}
\label{Le:independence}
For each $\alpha>0$, we have $\PPP{\mathcal{E}_n^c}=o(n^{-\alpha})$ as $n \to \infty$.
\end{lemma}

\begin{proof}
By subadditivity and independence, it follows that
\begin{align*}
\PPP{\mathcal{E}_n^c} & \leq \sum_{\jj\in \VV}\PPP{{\cal X}_n\cap \jj=\emptyset}\\
& = \sum_{\jj\in \VV}\big(  \PPP{X_1\not\in \jj} \big)^n\\
& = \sum_{\jj\in \VV} \big(1-\mu(\jj)\big)^n\\
& \leq \sum_{\jj\in \VV} \exp(-n\mu(\jj)).
\end{align*}
Here, the last inequality holds since $\log(1-x)\leq -x$ for each $x\in [0,1)$. Since $f\geq f_- >0$ on $K$, we have
$\mu(\jj)=\int_\jj f \text{d} \lambda \geq f_- \lambda(\jj)$, whence -- writing $\# M$ for the cardinality of a finite set $M$ --
\begin{align*}
\PP \big(\mathcal{E}_n^c\big) & \leq \sum_{\jj\in \VV} \exp\big(-n f_- \lambda(\jj)\big)\\
&\leq \#\VV \exp\big(- f_-(\log n)^{1+\varepsilon}\big).
\end{align*}
Since  $\#\VV\leq n/(\log n)^{1+\varepsilon}$, it follows that $n^\alpha \PPP{\mathcal{E}_n^c} \to 0  \text{ as }  n \to \infty$.
\end{proof}

\paragraph{The local property}
We now define a metric d on $\VV$ by putting d$(\jj,\jj'):=\max_{1\leq s\leq d}|j_{s}-j'_{s}|$ for any two different subcubes $\jj$ and $\jj'$, and
d$(\jj,\jj):= 0$, $\jj \in \VV$. Let $S(\jj, r) = \{\jj'\in \VV: \text{d}(\jj, \jj')\leq r\}$ be the ball of subcubes of radius $r$ centered at $\jj$.
For any $\jj\in \VV$, put
\[
M_\jj := \max_{i\leq n, X_i\in \jj} \mu(B(X_i,R_{i,n,k})),
\]
with the convention $M_\jj=0$ if ${\cal X}_n\cap \jj=\emptyset$. Conditionally on the event $\mathcal{E}_n$, and provided that d$(\jj,\jj')\geq 2k+1$,
the random variables $M_\jj$ and $M_{\jj'}$ are independent. Lemma \ref{Le:independence} is referred to as the \textit{asymptotic independence property}: conditionally on the event $\mathcal{E}_n$, which occurs with high probability, the extremes $M_\jj$ and $M_\jj'$ attained on two subcubes which are sufficiently distant from each other are independent.

The following lemma claims that, with high probability, two exceedances cannot occur in the same neighborhood.

\begin{lemma}
\label{Le:localproperty}
With the notation $a \wedge b := \min(a,b)$ for $a,b \in \mathbb{R}$, let
\[
R(n)= \sup_{\jj\in\VV}\sum_{i\neq i'\leq n} \PP \big( X_i,X_{i'}\in S(\jj,2k); \mu(B(X_i,R_{i,n,k})) \wedge \mu(B(X_{i'},R_{i',n,k}))>v_{n,k} \big).
\]
Then $R(n) = O(n^{-1}(\log n)^{2-d+\varepsilon})$ as $n \to \infty$.
\end{lemma}
Here, with a slight abuse of notation, we have identified the family of subcubes $S(\jj,2k)=\{\jj'\in \VV: \text{d}(\jj,\jj')\leq 2k\}$
with the set $\bigcup \big\{\jj': \jj'\in \VV \text{ and }  \text{d}(\jj,\jj')\leq 2k\}.$

We prepare the proof of Lemma \ref{Le:localproperty} with the following {result} that gives the volume of two $d$-dimensional balls.

\begin{lemma} \label{Le:balls}
	{If $x \in B(0,2)$ then}
	\begin{align*}
		\lambda(B(0,1)\cup B(x,1))=2\left(\kappa_d \left(1-\frac{\arccos(\|x\|/2)}{\pi}\right)+\frac{\|x\|\kappa_{d-1}}{2d}\left(\sqrt{1-(\|x\|/2)^2}\right)^{d-1}\right).
	\end{align*}
\end{lemma}

\begin{proof}
	We calculate the volume of $\lambda(B(0,1)\cup B(x,1))$ as the sum of the volumes of the following two congruent sets. The first one, say $B$, is given by the set of all points in $B(0,1)\cup B(x,1)$ that are closer to $0$ than to $x$ and for the second one we change the roles of $0$ and $x$. The set $B$ is the union of a cone $C$ with radius $\sqrt{1-(\|x\|/2)^2}$, height $\|x\|/2$ and apex at the origin and a set $D:=B(0,1)\setminus S$, where $S$ is a simplicial cone with external angle $\arccos(\|x\|/2)$. From elementary geometry, we obtain that the volumes of $C$ and $D$ are given by
	\begin{align*}
		\lambda(C)=\frac{\|x\|\kappa_{d-1}}{2d}\left(\sqrt{1-(\|x\|/2)^2}\right)^{d-1},\quad
		\lambda(D)=\kappa_d \left(1-\frac{\arccos(\|x\|/2)}{\pi}\right).
	\end{align*}
	This finishes the proof of the lemma.
\end{proof}

\begin{proof}[Proof of Lemma \ref{Le:localproperty}] For $z \in {[0,1]^d}$, let
	\begin{align*}
		r_{n,k}(z):=\inf\{r>0:\,\mu(B(z,r))> v_{n,k}\}.
	\end{align*}
	Writing $\#{\cal Y}(A)$ for the number of points of a finite set ${\cal Y}$ of random points in
	$\mathbb{R}^d$ that fall into a Borel set $A$, we have
	\begin{align*}
		\mu(B(z,R_{n,k}(z)))>v_{n,k}\quad \Longleftrightarrow \quad \# \mathcal{X}_n(B(z,r_{n,k}(z))) \le k-1.
	\end{align*}
	In the following, we assume that $r_{n,k}(X_{i'})\le r_{n,k}(X_{i})$ (which is at the cost of a factor $2$) and distinguish the two cases $X_{i'} \in B(X_i,r_{n,k}(X_{i}))$ and $X_{i'} \in S(\jj,2k) \setminus B(X_i,r_{n,k}(X_{i}))$. This distinction of cases gives
	\begin{align}
		&\PP \big( X_i,X_{i'}\in S(\jj,2k);  \mu(B(X_i,R_{i,n,k})) \wedge \mu(B(X_{i'},R_{i',n,k}))>v_{n,k} \big)\nonumber\\
		&\quad\le 2 \PP \big( X_i,X_{i'}\in S(\jj,2k); r_{n,k}(X_{i'})\le r_{n,k}(X_{i}); \mu(B(X_i,R_{i,n,k})) \wedge \mu(B(X_{i'},R_{i',n,k}))>v_{n,k} \big)\nonumber\\
		&\quad \le 2 \PP \big( X_i \in S(\jj,2k); X_{i'} \in  B(X_i,r_{n,k}(X_{i})); r_{n,k}(X_{i'})\le r_{n,k}(X_{i});\nonumber\\ &\quad \quad \quad \quad \mu(B(X_i,R_{i,n,k})) \wedge \mu(B(X_{i'},R_{i',n,k}))>v_{n,k} \big)\label{lemloccl}\\
		&\quad \quad +2 \PP \big( X_i\in S(\jj,2k),X_{i'} \in S(\jj,2k) \setminus B(X_i,r_{n,k}(X_{i})); r_{n,k}(X_{i'})\le r_{n,k}(X_{i});\nonumber\\
		&\quad \quad \quad \quad  \mu(B(X_i,R_{i,n,k}))\wedge \mu(B(X_{i'},R_{i',n,k}))>v_{n,k} \big).\label{lemlocfa}
	\end{align}
	We bound the summands \eqref{lemloccl} and \eqref{lemlocfa} separately. Since $X_i$ and $X_{i'}$  are independent, \eqref{lemloccl} takes the form
	\begin{align}
		&2\int_{S(\jj,2k)} \int_{B(x,r_{n,k}(x))}\mathds{1}\{r_{n,k}(y) \le r_{n,k}(x)\} \, \PP(\#(\mathcal{X}_n \setminus \{X_i,X_{i'}\} \cup \{x\})(B(y,r_{n,k}(y)))\le k-1; \nonumber\\
		&\quad \quad \#(\mathcal{X}_n \setminus \{X_i,X_{i'}\} \cup \{y\})(B(x,r_{n,k}(x)))\le k-1)\,\mu(\mathrm{d}y)\,\mu(\mathrm{d}x).\label{lemlocclint}
	\end{align}
	For $y \in B(x,r_{n,k}(x))$, the probability in the integrand figuring above is bounded from above by
	\begin{align}
		&\PP(\#(\mathcal{X}_n \setminus \{X_i,X_{i'}\})(B(y,r_{n,k}(y)))\le k-1;\nonumber\\
		&\quad  \#(\mathcal{X}_n \setminus \{X_i,X_{i'}\})(B(x,r_{n,k}(x)))\le k-2)\nonumber\\
		&\, \le \PP(\#(\mathcal{X}_n \setminus \{X_i,X_{i'}\})(B(y,r_{n,k}(y)))\le k-1;\nonumber\\
		&\, \quad  \#(\mathcal{X}_n \setminus \{X_i,X_{i'}\})(B(x,r_{n,k}(x))\setminus B(y,r_{n,k}(y))\le k-2).\label{negdep}
	\end{align}
	Since the random vector
	\[
	\left(\#(\mathcal{X}_n \setminus \{X_i,X_{i'}\})(B(y,r_{n,k}(y))),\#(\mathcal{X}_n \setminus \{X_i,X_{i'}\})(B(x,r_{n,k}(x))\setminus B(y,r_{n,k}(y)))\right)
	\]
	is negatively quadrant dependent (see \cite[Section 3.1]{JP83}), Equation \eqref{negdep} has the upper bound
	\begin{align}
		& \PP(\#(\mathcal{X}_n \setminus \{X_i,X_{i'}\})(B(y,r_{n,k}(y)))\le k-1)\nonumber\\
		&\quad \times \PP( \#(\mathcal{X}_n \setminus \{X_i,X_{i'}\})(B(x,r_{n,k}(x))\setminus B(y,r_{n,k}(y)))\le k-2)\nonumber\\
		&\, \le \PP(\#(\mathcal{X}_n \setminus \{X_i,X_{i'}\})(B(y,r_{n,k}(y)))\le k-1)\nonumber\\
		&\, \quad \times  \PP( \#(\mathcal{X}_n \setminus \{X_i,X_{i'}\})(B(x,r_{n,k}(x))\setminus B(y,r_{n,k}(x)))\le k-2),\label{ballest}
	\end{align}
	where the last inequality holds since $r_{n,k}(y) \le r_{n,k}(x)$. Analogously to Remark \ref{rem-byproduct},  the first probability {is}
	\begin{align*}
		\PP(\#(\mathcal{X}_n \setminus \{X_i,X_{i'}\})(B(y,r_{n,k}(y)))\le k-1)= \sum_{j=0}^{k-1} {n-2 \choose j} v_{n,k}^j (1-v_{n,k})^{n-2-j} \sim \frac{{\rm e}^{-t}}{n}.
	\end{align*}
	The latter probability in \eqref{ballest} is given by
	\begin{align}
		&\sum_{\ell=0}^{k-2} \binom{n-2}{\ell} \mu \big(B(x,r_{n,k}(x))\setminus B(y,r_{n,k}(x))\big)^\ell \big(1-\mu(B(x,r_{n,k}(x))\setminus B(y,r_{n,k}(x)))\big)^{n-2-\ell}. \label{bindist}
	\end{align}
	In a next step, we estimate $\mu \big(B(x,r_{n,k}(x))\setminus B(y,r_{n,k}(x))\big)$. Since $f(x) \ge f_- >0,\,x \in {[0,1]^d}$, and by the homogeneity of $d$-dimensional Lebesgue measure $\lambda$,
	we obtain
	\begin{align*}
		\mu \big(B(x,r_{n,k}(x))\setminus B(y,r_{n,k}(x))\big)&\ge f_- \lambda (B(x,r_{n,k}(x))\setminus B(y,r_{n,k}(x))) \\
		\quad& =  f_- r_{n,k}(x)^d  \lambda (B(0,1)\setminus B(r_{n,k}(x)^{-1}(y-x),1))\\
		\quad&=  f_- r_{n,k}(x)^d  \left(\lambda (B(0,1)\cup B(r_{n,k}(x)^{-1}(y-x),1))-\kappa_d\right).
	\end{align*}
	For $y \in B(x,r_{n,k}(x))$, Lemma \ref{Le:balls} {yields}
	\begin{align*}
		&\mu \big(B(x,r_{n,k}(x))\setminus B(y,r_{n,k}(x))\big)\nonumber\\
		&\quad \ge f_-  r_{n,k}(x)^{d} \left(\kappa_d\left(1-\frac{2\arccos(\|x-y\|/2r_{n,k}(x))}{\pi}\right)+\frac{\|x-y\|\kappa_{d-1}}{2dr_{n,k}(x)}\left(\sqrt{1-(\|x-y\|/2r_{n,k}(x))^2}\right)^{d-1}\right). \label{geoest}
	\end{align*}
	Since $\inf_{s>0} s^{-1}(1-2\arccos(s)/\pi)>0$, there is $c_0>0$ such that
	\begin{align*}
		\mu \big(B(x,r_{n,k}(x))\setminus B(y,r_{n,k}(x))\big) \ge c_0 \|x-y\| r_{n,k}(x)^{d-1} ,\quad x \in S(\jj,2k),\,y \in B(x,r_{n,k}(x)).
	\end{align*}
	{Equation} \eqref{bindist} and the bound $f(x) \le f_+,\,x \in {[0,1]^d},$ {give}
	\begin{align*}
		&\int_{B(x,r_{n,k}(x))}\mathds{1}\{r_{n,k}(y) \le r_{n,k}(x)\} \PP \big(\#(\mathcal{X}_n \setminus \{X_i,X_{i'}\})(B(x,r_{n,k}(x))\setminus B(y,r_{n,k}({x})))\le k-1\big)\, \mu(\mathrm{d}y)\\
		&\quad \le f_+ \sum_{\ell=0}^{k-2} \binom{n-2}{\ell} \int_{B(x,r_{n,k}(x))} \Big(c_0\|x-y\|r_{n,k}(x)^{d-1}\Big)^\ell \Big(1-c_0\|x-y\|r_{n,k}(x)^{d-1}\Big)^{n-2-\ell} \lambda(\mathrm{d}y).
	\end{align*}
	We now introduce spherical coordinates and obtain
	\begin{align*}
		&f_+ d \kappa_d \sum_{\ell=0}^{k-2} \binom{n-2}{\ell} \int_0^{r_{n,k}(x)} \Big(c_0 t r_{n,k}(x)^{d-1}\Big)^\ell \Big(1-c_0 t r_{n,k}(x)^{d-1}\Big)^{n-2-\ell} t^{d-1} \mathrm{d}t\\
		&\quad =f_+ d \kappa_d \sum_{\ell=0}^{k-2} \binom{n-2}{\ell} \int_0^{r_{n,k}(x)} \Big(c_0 t r_{n,k}(x)^{d-1}\Big)^\ell \exp \Big( (n-2-\ell) \log(1-c_0 t r_{n,k}(x)^{d-1})\Big) t^{d-1} \mathrm{d}t\\
		&\quad \le f_+d \kappa_d \sum_{\ell=0}^{k-2} \binom{n-2}{\ell} \int_0^{r_{n,k}(x)} \Big(c_0 t r_{n,k}(x)^{d-1}\Big)^\ell \exp \Big( -c_0(n-2-\ell)  t r_{n,k}(x)^{d-1}\Big) t^{d-1} \mathrm{d}t.
	\end{align*}
	Here, the last line follows from the inequality $\log s \le s-1,\,s>0$. Next,  we apply the change of variables
	\[
	t:=(c_0(n-2-\ell))^{-1} r_{n,k}(x)^{1-d} s \quad \Big(\text{i.e.}, s=c_0(n-2-\ell) t r_{n,k}(x)^{d-1}\Big),
	\]
	which shows that the last upper bound takes the form
	\begin{equation}\label{glaftersc}
		f_+d \kappa_d {c_0^{-d}} r_{n,k}(x)^{d(1-d)} \sum_{\ell=0}^{k-2} \binom{n-2}{\ell} (n-2-\ell)^{-d-\ell} \int_0^{c_0 (n-2-\ell)r_{n,k}(x)^d} s^{\ell+d-1} e^{-s} \,\mathrm{d}s.
	\end{equation}
	We now use the bounds $f_- \kappa_d r_{n,k}(x)^d \le v_{n,k}$, $\binom{n-2}{\ell}\le n^{\ell}/\ell!$, and the fact that the integral figuring in
	\eqref{glaftersc} converges as $n \to \infty$. Hence, the expression in \eqref{glaftersc} is bounded from above by  $c_1 n^{-1} (\log n)^{1-d}$, where $c_1$ is some
	positive constant. Consequently,  \eqref{lemloccl} is bounded from above by
	\begin{align}
		&c_1 n^{-1}(\log n)^{1-d} \lambda(S(\jj,2k)) \sup_{y \in S(\jj,2k)} \PP(\#(\mathcal{X}_n \setminus \{X_i,X_{i'}\})(B(y,r_{n,k}(y)))\le k-1)\nonumber\\
		&\quad  \sim c_2 n^{-3} (\log n)^{2-d+\varepsilon} \label{lemlocclbou}
	\end{align}
 for some $c_2>0$.
	
	By analogy with the reasoning above, \eqref{lemlocfa} is given by the integral
	\begin{align}
		&2\int_{S(\jj,2k)}\int_{S(\jj,2k) \setminus B(x,r_{n,k}(x))} \mathds{1}\{r_{n,k}(y) \le r_{n,k}(x)\} \,\PP\big(\#(\mathcal{X}_n \setminus \{X_i,X_{i'}\} \cup \{x\})(B(y,r_{n,k}(y))) \le k-1\big)\nonumber\\
		&\quad \times  \PP(\#(\mathcal{X}_n \setminus \{X_i,X_{i'}\})(B(x,r_{n,k}(x))\setminus B(y,r_{n,k}(y)))\le k-1)\,  \mu(\mathrm{d}y)\,\mu(\mathrm{d}x). \label{lemlocfaint}
	\end{align}
	If $y \notin B(x,r_{n,k}(x))$ and $r_{n,k}(x)\ge r_{n,k}(y)$, we have the lower bound
	\[
	\lambda(B(x,r_{n,k}(x))\setminus B(y,r_{n,k}(y))) \ge \frac{\lambda(B(x,r_{n,k}(x)))}{2}.
	\]
	Since {$f_+ \kappa_d r_{n,k}(x)^d \geq v_{n,k}$}, we find a constant $c_3>0$ such that
	\[
	\lambda(B(x,r_{n,k}(x))\setminus B(y,r_{n,k}(y))) \ge c_3 v_{n,k},
	\]
	whence
	\begin{align*}
		&\PP\big(\#(\mathcal{X}_n \setminus \{X_i,X_{i'}\})(B(x,r_{n,k}(x))\setminus B(y,r_{n,k}(y)))\le k-1\big)\\
		&\quad \le \sum_{\ell=0}^{k-1} \binom{n-2}{\ell} \big(c_3v_{n,k}\big)^\ell \big(1-c_3v_{n,k}\big)^{n-2-\ell}\\
		&\quad \sim \frac{c_3^{k-1}}{(k-1)!} \big(\log n\big)^{k-1} \exp \big(n \log (1-c_3v_{n,k})\big)
	\end{align*}
	as $n \to \infty$. Since $\log s \le s-1$ for $s>0$, \eqref{lemlocfaint} is bounded from above by
	\begin{align}
		&c_4 n^{-c_3} \lambda(S(\jj,2k))^2 \sup_{y \in S(\jj,2k)} \PP(\#(\mathcal{X}_n \setminus \{X_i,X_{i'}\})(B(y,r_{n,k}(y)))\le k-1)\nonumber\\
		&\quad \sim c_5 (4k+1)^{2d} \frac{(\log n)^{2+2\varepsilon}}{n^{3+c_3}}, \label{lemlocfabou}
	\end{align}
where $c_4$ and $c_5$ are positive constants.  {Summing over all $i\neq i'\leq n$, it follows} from  \eqref{lemlocclbou} and  \eqref{lemlocfabou} 
	that $R(n)=O(n^{-1}(\log n)^{2-d+\varepsilon})$ as $n \to \infty$, {which} finishes the proof of Lemma~\ref{Le:localproperty}.
\end{proof}

\paragraph{ A Poisson approximation result based on the Chen-Stein method}
In this paragraph, we recall a Poisson approximation result due to Arratia \textit{et al.} \cite{AGG}, which is based on the Chen-Stein method.
To this end, we consider a finite or countable collection $(Y_\alpha)_{\alpha\in I}$ of $\{0,1\}$-valued random variables and we let $p_\alpha = \PP(Y_\alpha=1)>0$, $p_{\alpha\beta}=\PP(Y_\alpha=1, Y_\beta=1)$. Moreover, suppose that
for each $\alpha\in I$, there is a set $B_\alpha\subset I$ that contains $\alpha$. The set $B_\alpha$ is regarded as a neighborhood
of $\alpha$ that consists of the set of indices $\beta$ such that $Y_\alpha$ and $Y_\beta$ are \textit{not} independent. Finally, put
\begin{equation}
\label{eq:defb}
b_1=\sum_{\alpha\in I} \sum_{\beta \in B_\alpha}p_\alpha p_\beta, \quad b_2=\sum_{\alpha\in I}\sum_{\alpha\neq \beta \in B_\alpha}p_{\alpha\beta},
\quad b_3=\sum_{\alpha\in I}\E \big[\left| \E[ Y_\alpha-p_\alpha|\sigma(Y_\beta: \beta\not\in B_\alpha)]  \right| \big] .
\end{equation}

\begin{theorem} {\rm (}Theorem 1 of \cite{AGG} {\rm )}
\label{Th:AGG}
Let $W=\sum_{\alpha \in I}Y_\alpha$, and assume $\lambda: =\E(W) \in (0,\infty)$. Then
\[
d_{TV}(W, {{\rm{Po}}}(\lambda)) \leq 2(b_1+b_2+b_3).
\]
\end{theorem}

\paragraph{Proof of Theorem \ref{Th:poissonapprox}}
Recall $v_{n,k}$ from \eqref{defynk} and $C_{n,k}$ from \eqref{defcnk}.  Put
\[
\widehat{C}_{n,k} := \sum_{\jj\in \VV} \ind \big\{ M_\jj>v_{n,k} \big\}.
\]

The following lemma claims that the number {$C_{n,k}$} of exceedances  is close to the number of subcubes for which there exists at least one exceedance, i.e. $\widehat{C}_{n,k}$, and that $\widehat{C}_{n,k}$ can be approximated by a Poisson random variable.

\begin{lemma} \label{Le:mainlemma} We have
\begin{enumerate}
\item[{{\rm a)}}] $\PP (C_{n,k} \neq \widehat{C}_{n,k} ) =  O\left(\left(\log n\right)^{1-d}   \right)$,
\item[{{\rm b)}}] $d_{TV}(\widehat{C}_{n,k}, {{\rm{Po}}}(\E[\widehat{C}_{n,k}]))\ =  O\left(\left(\log n\right)^{1-d}   \right)$,
\item[{{\rm c)}}]  $\E [ \widehat{C}_{n,k}] = {\rm e}^{-t} +   O\left( \frac{\log\log n}{\log n}  \right)$.
\end{enumerate}
\end{lemma}


\begin{proof}
Assertion a) is a direct consequence of Lemma~\ref{Le:localproperty}  and of the inequalities
\begin{align*}
\PP(C_{n,k}\neq \widehat{C}_{n,k})  & = \PP \big(\exists \jj\in\VV, \exists i,\ell \text{ s.t. } X_i,X_\ell \in \jj;\; \mu(B(X_i,R_{i,n,k})) \wedge \mu(B(X_{\ell},R_{\ell,n,k}))>v_{n,k} \big)\\
& \leq \sum_{\jj\in \VV}\sum_{i\neq \ell \leq n}\PPP{X_i,X_{\ell}\in \jj;\; \mu(B(X_i,R_{i,n,k})) \wedge \mu(B(X_{\ell},R_{\ell,n,k}))>v_{n,k} }\\
& \leq \frac{n}{(\log n)^{1+\varepsilon}}\times R(n).
\end{align*}

\bigskip

To prove b), we apply Theorem~\ref{Th:AGG} to the collection $(Y_\alpha)_{\alpha\in I} = (M_{\jj})_{\jj\in \VV}$. Recall that, conditionally on the event $\mathcal{E}_n$,
the random variables $M_{\jj}$ and $M_{\jj'}$ are independent provided that d$(\jj,\jj')\geq 2k+1$. With a slight abuse of notation, we omit to condition on $\mathcal{E}_n$ since this event
occurs with probability tending to 1 as $n \to \infty$ (Lemma \ref{Le:independence}) at a rate which is at least polynomial. The first two terms in \eqref{eq:defb} are
\[
b_1=\sum_{\jj\in \VV} \sum_{\jj' \in S(\jj,2k)}p_\jj p_{\jj'}, \quad b_2=\sum_{\jj\in \VV}\sum_{\jj\neq \jj' \in S(\jj,2k)}p_{\jj\jj'}, \]
where
\[
p_{\jj}= \PP(M_\jj > v_{n,k}), \quad p_{\jj\jj'}=\PP(M_\jj>v_{n,k}, M_{\jj'}>v_{n,k}).
\]
The term $b_3$ figuring in \eqref{eq:defb} equals 0 since, conditionally on $\mathcal{E}_n$, the random variable $M_\jj$ is independent of the $\sigma$-field
$\sigma(M_{\jj'}: \jj'\not\in S(\jj,2k))$.  Thus, according to Theorem~\ref{Th:AGG}, we have
\[d_{TV}(\widehat{C}_{n,k}, \text{Po}({\E[\widehat{C}_{n,k}]}))\leq 2(b_1+b_2).\]
{Fist, we deal with $b_1$}. As for the first assertion, notice that for each $\jj\in \VV$, using symmetry, we obtain
\begin{align*}
p_\jj & = \PP \bigg(\bigcup_{i\leq n}\left\{X_i\in \jj, \mu(B(X_i,R_{i,n,k})) >v_{n,k}    \right\}\bigg) \\
& \leq n \cdot \PP \big( X_1\in \jj, \mu(B(X_1,R_{1,n,k})) >v_{n,k} \big)\\
& = n \cdot \int_{\jj} \PP(\mu(B(x,R_{1,n,k})) >v_{n,k}|X_1=x) f(x) \, \text{d} x\\
& \le n f^+ \lambda(\jj) \int_{\jj} \PP(\mu(B(x,R_{1,n,k})) >v_{n,k}|X_1=x)  \, \frac{1}{\lambda(\jj)} \, \text{d} x\\
& = n f^+ \lambda(\jj) \PP(\mu(B(\widetilde{X}_1,R_{1,n,k})) >v_{n,k}),
\end{align*}
where $\widetilde{X}_1$ is independent of $X_2,\ldots,X_n$ and has a uniform distribution   over $\jj$.
Invoking Remark~\ref{rem-byproduct}, the probability figuring in the last line is asymptotically equal to ${\rm e}^{-t}/n$ as $n \to \infty$.
Since  $\lambda(\jj)=O\left( (\log n)^{1+\varepsilon}/n\right)$, we thus have
\[
p_\jj \le C \cdot \frac{(\log n)^{1+ \varepsilon}}{n},
\]
where $C$ is a constant that does not depend on $\jj$.
 Since $\#\VV\leq \frac{n}{ (\log n)^{1+\varepsilon}}$ and $\#S(\jj, 2k)\leq (4k+1)^d$, summing over $\jj,\jj'$ gives
\[
b_1  \leq  C^2 \sum_{\jj\in \VV}\sum_{\jj'\in S(\jj,2k)}   \left( \frac{(\log n)^{1+\varepsilon}}{n}  \right)^2
 = O\left(  \frac{(\log n)^{1+\varepsilon}}{n} \right).
\]

\medskip

To {deal with $b_2$}, notice that, for each $\jj,\jj'\in \VV$ and $\jj'\in S(\jj,2k)$, we have
\begin{align*}
p_{\jj\jj'} & = \PP \Big( \bigcup_{i\neq i'\leq n}  \left\{  X_i\in \jj, X_{i'}\in S(\jj,2k),  \mu(B(X_i,R_{i,n,k}))\wedge \mu(B(X_{i'},R_{i',n,k})) >v_{n,k}  \right\}   \Big)\\
& \leq \PP \Big( \bigcup_{i\neq i'\leq n}  \left\{  X_i, X_{i'}\in S(\jj,2k); \mu(B(X_i,R_{i,n,k})) \wedge  \mu(B(X_{i'},R_{i',n,k})) >v_{n,k}  \right\}    \Big).
\end{align*}
Using subadditivity, and taking the supremum, we obtain
\begin{align*}
b_2 & \leq \sum_{\jj\in \VV}\sum_{\jj'\in S(\jj,2k)} \sup_{\jj\in\VV}\sum_{i\neq i'\leq n}\PP \big(X_i,X_{i'}\in S(\jj,2k), \mu(B(X_i,R_{i,n,k})) \wedge\mu(B(X_{i'},R_{i',n,k}))>v_{n,k}\big) \\
& \leq \frac{n}{(\log n)^{1+\varepsilon}}\times (4k+1)^d\times R(n).
\end{align*}
According to  Lemma \ref{Le:localproperty}, the last term equals  $O\left(\left(\log n\right)^{1-d}   \right)$, {which concludes the proof of b).}

\bigskip

To prove c), observe that
\[
\big{|} \E [\widehat{C}_{n,k}] - {\rm e}^{-t} \big{|} \le \big{|}\E[\widehat{C}_{n,k}] -  \E[C_{n,k}]\big{|} + \left|\E[C_{n,k}] - {\rm e}^{-t}\right|.
\]
By Theorem~\ref{thmmetricspace}, the last summand is $O\left( \frac{\log\log n}{\log n}  \right)$. Since $C_{n,k} \ge \widehat{C}_{n,k}$, we further have
\begin{align*}
\big{|}\E [\widehat{C}_{n,k}] - \E[C_{n,k}]\big{|} & = \E [C_{n,k}-\widehat{C}_{n,k}]\\
& = \E \bigg( \sum_{i\leq n} \ind \{\mu(B(X_i, R_{i,n,k}))>v_{n,k}\} - \sum_{\jj\in\VV} \ind \{M_\jj>v_{n,k}\}  \bigg)\\
& = \sum_{\jj\in\VV} \E \bigg[ \Big( \sum_{i\leq n} \ind\{X_i\in \jj\} \ind \{\mu(B(X_i,R_{i,n,k})) > v_{n,k} \}  - 1 \Big)\ind \{M_\jj>v_{n,k}\} \bigg]\\
& \leq \sum_{\jj\in \VV}\sum_{i\neq i'\leq n} \mathbb{P} \big(X_i,X_{i'}\in \jj, \mu(B(X_i,R_{i,n,k})), \mu(B(X_{i'},R_{i',n,k}))> v_{n,k} \big)\\
& \leq \#\VV \times R(n).
\end{align*}
According to  Lemma \ref{Le:localproperty}, the last term equals  $O\left(\left(\log n\right)^{1-d}   \right)$. This concludes the proof of Lemma \ref{Le:mainlemma} and thus of Theorem~\ref{Th:poissonapprox}.

\end{proof}

\section{Concluding remarks}
When dealing with limit laws for large $k$th-nearest neighbor {\em distances} of a sequence of i.i.d. random points in $\mathbb{R}^d$ with density $f$, which take values in a bounded region $K$,
the modification of the $k$th-nearest neighbor distances made in \eqref{gumbelknnbound} (by introducing the ``boundary distances'' $\|X_i - \partial K\|$)
and  the condition that $f$ is bounded away from zero, which have been adopted in \cite{HE82} and \cite{HE83}, seem to be crucial, since boundary effects play
a decisive role (\cite{DH89, DH90}). Regarding $k$th-nearest neighbor balls with {\em large probability volume}, there is no need to introduce $\|X_i - \partial K\|$.
It is an open problem, however, whether Theorem~\ref{Th:poissonapprox} continues to hold for densities that are not bounded away from zero.

A second open problem refers to Theorem~\ref{thmmetricspace}, which states convergence of expectations of $C_{n,k}$ in a setting beyond the finite-dimensional case.
Since $C_{n,k}$ is non-negative, the sequence $(C_{n,k})_{k}$ is tight by Markov's inequality. Can one find conditions on the underlying distribution that ensure
convergence in distribution to some random element of the metric space?

\end{document}